\documentclass[11pt,a4paper]{article}
\usepackage{amsfonts}
\usepackage{amsmath}
\usepackage{amssymb}

\def\eqref#1{(\ref{#1})}
\newcounter{theoremcount}
\newcommand{\thrm}{\par\refstepcounter{theoremcount}%
{\bf Theorem \arabic{theoremcount}.} }

\newcounter{lemmacount}
\newcommand{\lemm}{\par\refstepcounter{lemmacount}%
{\bf Lemma \arabic{lemmacount}.} }

\newcommand{\Rd}{{\mathbb R}^d}

\newcommand{\R}{{\mathbb R}}

\newcommand{\N}{{\mathbb N}}
\newcommand{\proof}{{\bf Proof. } }

\newcommand{\B}{\mathcal B}
\newcommand{\E}{\mathsf E}
\newcommand{\lm}{\lambda}
\newcommand{\sgm}{\sigma}
\begin{document}
\begin{center}
{\Large
Existence And Uniqueness Of Stationary Solution Of
Nonlinear Stochastic Differential Equation With Memory.
\\}

\bigskip

\bf Yuri Bakhtin\footnote{ International Institute
                   of Earthquake  Prediction Theory and
                   Mathematical Geophysics, Moscow, e-mail: bahtin@mitp.ru, bakhtine@mail.ru }
\end{center}
\section{Introduction.}

In this paper a stochastic differential equation (SDE) with infinite memory is considered.
The drift coefficient of the equation is a nonlinear functional of the past history of the
solution. Sufficient conditions for existence and uniqueness of stationary solution are given.
This work is motivated by recent papers
 \cite{ems} and  \cite{el} where stochastically forced nonlinear equations of hydrodynamics
were considered and it was shown how the infinite-dimensional stochastic  Markovian
dynamics related to these equations can be reduced to finite-dimensional stochastic dynamics.
The corresponding  finite-dimensional systems are however essentially non-Markovian.

So, the important problem of existence and uniqueness of stationary solutions for
stochastic hydrodynamical equations is tightly related to existence and uniqueness
of stationary solutions of SDEs with infinite memory.
Some results in the area were established in \cite{itonisio}. In the first part of this
paper some necessary notions are introduced and the main result is stated.
A proof of the main result is given in the second part. We combine the approach of \cite{itonisio}
with an interesting method for establishing the desired uniqueness suggested in \cite{ems} and
\cite{el} for the problems considered therein.

The equation under consideration is
\begin{equation}
dX(t) = a(\pi_tX)dt + dW(t).
\label{sdde}
\end{equation}

Here $W(t),t\in\R$ is standard $d-$dimensional Wiener process
(i.e. a Gaussian $\Rd$-valued stochastic process
with continuous trajectories
defined on the whole real line $\R$ with
independent and stationary increments, $W(0)=0$, $\E W(t) = 0$,
$\E W(t)^2=|t|, t\in\R$),
$\pi_t$ is a map from the space
$C$ of $\Rd$-valued continuos functions defined on $\R$ to the space
 $C_-$ of continuous functions defined on $\R_-=(-\infty,0]$:
$$\pi_tX(s)=X(s+t),\quad s\in\R_-.$$
This map gives the past history of a continuous process up to time $t\in\R$.
From now on suppose $a(\cdot):C_-\to\Rd$ to be a continuous functional with respect
to metric
$$
\rho_-(f,g)=\sum_{n=1}^{\infty}2^{-n}(\|f-g\|_n\wedge1),\quad f,g\in C_-,
$$
which defines LU-topology on the space $C_-$. Here
  $\|h\|_n=\max_{-n\le t\le 0}|h(t)|,$ and  $|\cdot|$ denotes the Euclidean norm.

For a stochastic process $X$ and a set $A\subset\R$ the
$\sgm$-algebra generated by r.v.'s $X(s), s\in A$ will be denoted by
$\sgm_A(X)$ and the $\sgm$-algebra generated by r.v.'s
$X(s)-X(t), s,t\in A$ will be denoted by
$\sgm_A(dX)$.

Consider the space $\Omega=C\times C$ with the metric analogous to
the metric $\rho_-$ defined above.

A probability measure $P$ on the space $\Omega$ with Borel
$\sgm$-algebra $\B$
is said to define a solution to the
equation \eqref{sdde} on $\R$ if the following three conditions are
fulfilled with respect to the measure $P$:

1. The projection $W:C\times C\to C$, $\omega = (\omega_1,\omega_2)\mapsto\omega_2$,
is a standard $d$-dimensional Wiener process.

2. For any $t\in\R$
\begin{equation}
\sgm_{(-\infty,t]}(X)\vee\sgm_{(-\infty,t]}(dW)\quad \mbox{is independent of}\quad
\sgm_{[t,\infty)}(dW).
\label{independence}
\end{equation}
Here and further $X:C\times C\to C$, $\omega = (\omega_1,\omega_2)\mapsto\omega_1$.

3. If $s<t$ then
\begin{equation}
X(t)-X(s)\stackrel{\mbox{\small a.s.}}{=}\int_s^ta(\pi_\theta X)d\theta+W(t)-W(s).
\label{fulfilsdde}
\end{equation}

If in addition the distribution of the process
$$
(X,dW) \equiv (X(t), -\infty<t<\infty, W(v)-W(u), -\infty<u<v<\infty)
$$
does not change under time shifts then the measure $P$ is said to define a
stationary solution.

Let's state the main result.

\thrm {\it
Let the drift coefficient $a(\cdot)$ satisfy  the following conditions:

{\rm 1.} There exist such constants $K>0,\lm>0$ that the estimate
\begin{equation}
|a(x_-)-a(y_-)|\le K  \int_{-\infty}^0e^{\lm t}|x_-(t)-y_-(t)|dt
\label{lipschits}
\end{equation}
is fulfilled whenever $x_-,y_-\in C_-$, $x_-(0)=y_-(0)$ and the integral in the right-hand side is finite.

{\rm 2.} There exist such constants $C_1\ge0$ and $C_2>0$ that
\begin{equation}
(a(x_-),x_-(0))\le C_1-C_2|x_-(0)|^2,\quad x_-\in C_-.
\label{khas'minski}
\end{equation}

{\rm 3.} There exist such a constant $C_3>0$ that
\begin{equation}
|a(x_-)| \le C_3|x_-(0)|,\quad x_-\in C_-.
\label{lineargrowth}
\end{equation}

Then there exist a probabilistic measure $P$ on the space $C\times C$
which defines a stationary solution of the equation {\rm \eqref{sdde}}.
Such measure is unique in the class of measures for which almost every
realization $X$ possesses the following property:
\begin{equation}
|X(t)|\le K'e^{\lambda'|t|},\quad t\le0.
\label{subexp}
\end{equation}
Here $K'\in \R$ and $\lambda'\in(0,\lambda)$ are some constants depending on the realization $X$.

}

\section{Proof of the main result.}

First, let us prove the existence of the stationary solution using
the Krylov--Bogolyubov approach.

A probabilistic law in
$C\times C$ is said to define a solution of Cauchy problem for the
equation \eqref{sdde} with initial data $x_-\in C_-$
if  the following conditions are satisfied:
 $W$ is a standard Wiener process;
for every $t\in\R$ the relation
\eqref{independence} is true;
the equality \eqref{fulfilsdde} is fulfilled for
 $s=0$ and every $t>0$; $X(t)=x_-(t)$ for any $t<0$. Existence theorem for solutions
 of Cauchy problem is proved in \cite{itonisio}.

Let $P_0$ denote such a law for the initial data identically equal to zero and
$P_s$ denote the time $s$-shift of this distribution i.e. a solution
of the Cauchy problem subject to zero initial data defined on the
set $(-\infty,s],s\in\R$. Formally $P_s=P_0\theta_s^{-1}$ where
$\theta_s(f,g) = (\widetilde f,\widetilde g)$,
$\widetilde f(t) = f(t-s)$, $\widetilde{ g}(t) = g(t-s)-g(-s)$.

Since the function $P_s(E)$ is  measurable with respect to $s$ for all $E\in\B$ (see \cite{itonisio}),
for $T>0$ one can define a propability measure
$$
Q_T(\cdot)=\frac{1}T\int_{-T}^0 P_s(\cdot)ds
$$
on the space $(\Omega,\B)$.
We will show that the family of measures $\{Q_T\}$ is tight.

Theorem 12.3 of the book \cite{bil} implies that
in order to prove tightness of a family of measures
$\{Q_T\}$ it is sufficient to verify that corresponding
one-dimensional distributions constitute a tight family and there exist
a nondecreasing continuous function and constants
 $\gamma\ge0, \alpha>1, z_0>0$ such that
for all
$z\in(0,z_0)$, $t_1,t_2\in\R$ and  $s\ge0$
the following estimate holds:
$$
Q_T\{|X(t_2)-X(t_1)|\ge z\}\le z^{-\gamma}|F(t_2)-F(t_1)|^\alpha.
$$
Introduce stopping time $\tau_r(t)=t\wedge \tau_r$ where
$\tau_r=\inf\{t:|X(t)|\ge r\}$.
Then for $T>0$ the Ito formula implies that the equality
$$
X^2(\tau_r(T))=2\int_0^{\tau_r(T)}(a(\pi_\theta(X),X(\theta))d\theta +
2\int_0^{\tau_r(T)}(X(\theta),dW(\theta)) + \tau_r(T)
$$
holds $P_0$-a.s.
Taking expectations of both sides, passing to
limit $r\to\infty$, using the regularity of the solution
($\tau_r\to\infty$ a.s. for $r\to\infty$) and
finiteness of second-order moments of the solution established
in \cite{itonisio} and the inequality \eqref{khas'minski}, one can obtain that
for all $S\in\R$
\begin{multline*}
\E_{P_0} |X(T)|^2\le (2C_1+1)T-2C_2\int_0^T\E_{P_0}|X(\theta)|^2d\theta
\\ \le (2C_1+1)T-2C_2\int_{S-T}^S\E_{P_\eta}|X(S)|^2d\eta.
\end{multline*}
Dividing both parts of this inequality by $T$ and considering the last summand, one
can obtain that for some positive constant $M$ and for all $S>0$
\begin{equation}
\limsup_{T\to\infty}\sup_{t\in[-S,S]}\E_{Q_T} |X(t)|^2\le M.
\label{bounded2ndmoment}
\end{equation}

Now let us estimate increments of the process $X$.
\begin{multline}
\label{increments}
Q_T\{|X(t_2)-X(t_1)|> z\}\\ \le  Q_T\{|W(t_2)-W(t_1)|>z/2\}+
Q_T\left\{\int_{t_1}^{t_2}a(\pi_{\theta}X)d\theta >z/2 \right\}
\\\le
16z^{-4}\E_{Q_T}|W(t_2)-W(t_1)|^4+4z^{-2}\E_{Q_T}\left(\int_{t_1}^{t_2}
a(\pi_{\theta}X)d\theta\right)^2.
\end{multline}
The next inequality is a consequence of the Fubini theorem, elementary
inequality $|xy|\le(x^2+y^2)/2$, well-known expression for moments of
Gaussian distribution and relations \eqref{lineargrowth} and \eqref{increments}:
\begin{equation}
Q_T\{|X(t_2)-X(t_1)|> z\}\le 48z^{-4}|t_2-t_1|^2 + 4C_3^2z^{-2}M|t_2-t_1|^2.
\label{incrementsfinal}
\end{equation}
Tightness of the family of projections of measures $Q_T$ on the
first component $\{Q_TX^{-1}\}$ and hence the desired tightness of the
family $\{Q_T\}$ is implied now by
\eqref{bounded2ndmoment} and \eqref{incrementsfinal}.

So,
$Q_{T_n}\stackrel{Law}{\to} Q_\infty$ when $n\to\infty$
for some sequence $(T_n)_{n\in\N}$ and results of
\cite{itonisio} imply that  $Q_\infty$ defines a stationary
solution of the equation \eqref{sdde}.

\lemm\label{bounds} {\it For any $\delta>0$ the following estimate is true $P-$a.s.
$$
\lim_{t\to-\infty}|X_t|\cdot|t|^{1/2+\delta}=0.
$$
}

\proof An estimate for measure $P$, analogous to the estimate \eqref{incrementsfinal},
imlies that for any  $s\in\R$
$$
P\left\{\max_{t\in[s,s+1]}|X(t)|>Kz\right\}\le P\{X(s)>z\}
+P\{X(s+1)>z\}+C(z^{-2}+z^{-4})
$$
for sufficiently large constants $K,C>0$.
Using this inequality and Chebyshev inequality and uniform in $t\in\R$
boundedness
of the second-order moment of  $X(t)$
one obtains that for all
$\delta_0\in(0,\delta)$ the series
$$
\sum_{n=0}^{\infty} P\left\{\max_{t\in[-n-1,-n]}
|X(t)|>Kn^{1/2+\delta_0}\right\}
$$
is convergent and the lemma follows from the Borel--Kantelly lemma.

In particular, Lemma \ref{bounds} implies that the trajectories of
the process $X$ satisfy the condition \eqref{subexp} $P$-a.s.

Now we turn to the proof of uniqueness. Consider an arbitrary measure
$P$ which defines a stationary solution of the equation \eqref{sdde}.
Suppose also that the realizations of the process $X$ satisfy condition
\eqref{subexp} $P$-a.s. Introduce a space $C_+$ of $\Rd$-valued
continuous functions
defined on $\R_+=[0,\infty)$. For
$x_-\in C_-$ we denote $P_{x_-}$ the measure on $\Omega_+ = C_+\times C_+$
which defines a solution of Cauchy problem with the initial data $x_-$.
$P_{x_-}$ is a conditional distribution of the measure $P$ conditioned on
$X_-=x_-$.

\lemm\label{samepastequiv}
{\it Condition {\rm 1} of Theorem {\rm 1} implies that there exists a set
$A\subset C_-$ such that $P(\pi_0X\in A)=1$ and if
$x_-,y_-\in A$ и $x_-(0) = y_-(0)$ then
the measures  $P_{x_-}$
and $P_{y_-}$ are equivalent.
}

\proof Consider $x_-, y_-\in C_-$ such that each of these functions
admits an exponential estimate like \eqref{subexp}.
To prove that
$P_{y_-}$ is absolutely continuous with respect to $P_{x_-}$,
we use the Girsanov theorem and verify th Novikov condition
(see, e.g., \cite[Chapter 8]{RevuzYor}).
The same reasoning will be valid for interchanged $x_-$ and $y_-$.

The Novikov condition can be written as follows:
\begin{equation}
\E_{P_{x_-}}\exp\left\{\frac12\int_0^\infty
\left|a(\pi_tX)-a(\pi_t(\widehat{yX}))\right|^2dt\right\}<\infty
\label{novikov}
\end{equation}
where
$$
\begin{array}{lcl}
\widehat{yX}(s)&=&\left\{\begin{array}{ll}
y_-(s),& s<0\\
X(s),& s\ge0.
\end{array}
\right.
\end{array}
$$
The condition \eqref{subexp} implies that for some constants $K'>0$ и $\lambda' \in(0,\lambda)$
the inequality $|x_-(s)-y_-(s)|\le K'e^{\lambda' |s|}$ is fulfilled, and
the condition \eqref{lipschits} implies that
\begin{multline*}
\left|a(\pi_tX)-a(\pi_t(\widehat{yX}))\right|\le
K\int_{\infty}^0|x_-(s)-y_-(s)|e^{\lambda(s-t)}ds
\\ \le KK'e^{-\lambda t}\int_{\infty}^0e^{(\lambda-\lambda') s}ds\le Le^{-\lambda t}
\end{multline*}
for some constant $L>0$.
So,  in \eqref{novikov} the expectation of a bounded random variable is taken.
The lemma is proved.

Lemma \ref{samepastequiv} implies the following result:

\lemm\label{CommonComponent} {\it For $P\{X(0)\in\cdot\}$-almost every $l\in\Rd$
the measure $P\{\pi_0^+X\in\cdot\ |X(0)=l\}$ has a component which is equivalent to
the measure
$P\{\pi_0^+X\in\cdot\ |\pi_0X\equiv l\}$. Here $\pi_t^+X(s)=X(s+t), s\in\R_+.$}

In the same way the following lemma can be proved.
\lemm {\it There exists a set $A_-\subset C_-$ such that
$P\{\pi_0 X \in A_-\}=1$ and if  $x_-\in A_-$ then the
projection of the measure $P_{x_-}$ on the space $\Omega_{[0,T]}$
of $\Rd\times\Rd$-valued continuous functions defined on finite
segment is equivalent to the distribution of the
standard Wiener process in $\R^{2d}$  emitted from $(x_-(0),0)$.}

Let $R_{T,x_-}(\cdot)$ be the restriction of this measure on sets
of the form $\{z\in\Omega_{[0,T]}| z(T)\in B\}$, $B\in\B(\R^{2d})$.
Since restrictions of equivalent measures on $\sigma$-subalgebra
are equivalent, the measure $R_{T,x_-}(\cdot)$ is equivalent to
a non-degenerate Gaussian measure in $\R^{2d}$, and hence to the Lebesgue measure.
This fact and stationarity of the process $X$ imply the following result.

\lemm \label{equivTOlebesgue}
{\it For any $t\in\R$ the measure $P\{X(t)\in\cdot\}$ has a component which is
equivalent to the Lebesgue measure.
}

The proof of uniquemess of stationary solution of the equation \eqref{sdde}
given here is based on the lemmas above. It is a modification of
reasoning from \cite{ems}.

Suppose there are two different ergodic measures
 $P^{(1)}$ и $P^{(2)}$ defining stationary solutions.
There exists a bounded functional $F$ such that
$$
\overline F_1 = \E_{P^{(1)}}
F((X,W))\ne\E_{P^{(2)}}F((X,W)) = \overline F_2,
$$
and for some $S>0$
  $x(s)=y(s), s\in[0,S]$ implies $F(x)=F(y)$.

Then there exist sets $B_1,B_2\in\B$ such that $P^{(i)}(B_i)=1$ and
$$
\lim_{T\to\infty}\frac1T\int_0^T F(\theta_s(X,W))ds = \overline F_i
 \mbox{ on } B_i,i=1,2.
$$
Lemmas  \ref{CommonComponent} and \ref{equivTOlebesgue} imply that
$$
P^{(2)}(B_1) = \int_{\Rd} P^{(2)}(B_1|X(0)=l)P^{(2)}(X(0)\in dl) > 0.
$$
So, $P^{(2)}(B_1\cap B_2) > 0$ and $B_1\cap B_2\ne\varnothing$, which
contradicts the assumption $\overline F_1\ne \overline F_2$.

The author is grateful to Professor Ya.G.Sinai for statement of the
problem and useful discussions.

\end{document}